\documentclass[12pt, final]{l4dc2022}

\title[Heterogenous Multi-agent Learning]{On the Heterogeneity of Independent Learning Dynamics in \\Zero-sum Stochastic Games}
\usepackage{times}
\usepackage{color}

\usepackage{algorithm}
\usepackage{algpseudocode}

\newtheorem{assumption}[theorem]{Assumption}

\DeclareMathOperator*{\argmax}{argmax}

\newcommand{\be}{\begin{equation}}
\newcommand{\ee}{\end{equation}}
\newcommand{\nn}{\nonumber}
\newcommand{\val}{\mathrm{val}}
\renewcommand{\tilde}{\widetilde}
\renewcommand{\hat}{}

\author{%
 \Name{Muhammed O. Sayin} \Email{sayin@ee.bilkent.edu.tr}\\
 \addr Bilkent University, Electrical and Electronics Engineering Department, Ankara, Turkey 06800
 \AND
 \Name{K. Alperen Cetiner} \Email{alperen.cetiner@bilkent.edu.tr}\\
 \addr Bilkent University, Electrical and Electronics Engineering Department, Ankara, Turkey 06800\\
 Aselsan A.\c{S}., Ankara, Turkey 06800
}

\begin{document}

\maketitle

\begin{abstract}%
We analyze the convergence properties of the two-timescale fictitious play combining the classical fictitious play with the $Q$-learning for two-player zero-sum stochastic games with player-dependent learning rates. We show its almost sure convergence under the standard assumptions in two-timescale stochastic approximation methods when the discount factor is less than the product of the ratios of player-dependent step sizes. To this end, we formulate a novel Lyapunov function formulation and present a one-sided asynchronous convergence result.
\end{abstract}


\section{Introduction}

Multi-agent reinforcement learning has become the frontier of many advancements in artificial intelligence systems, where autonomous agents make decisions in dynamic environments (e.g., see \cite{ref:Zhang19} and the references therein). Heterogeneity and independence of the learning dynamics adopted by these autonomous agents are inevitable in practical applications of multi-agent systems. However, there has been very limited progress addressing it in the context of stochastic games (introduced by \cite{ref:Shapley53}) - a canonical model for dynamic multi-agent interactions.

Independent learning in strategic-form games played repeatedly has been studied extensively with many well-established results, e.g., see \citep{ref:Fudenberg98,ref:Young04,ref:Fudenberg09}. On the other hand, for stochastic games, \cite{ref:Arslan17,ref:Wei17} presented learning dynamics with double-loop-like update rules necessitating coordination among players that may not be inline with their best interests. Recently, \cite{ref:Leslie20} has drawn a two-timescale learning framework in which continuous-time best-response dynamics could also convergence to an equilibrium of a two-player zero-sum stochastic game though the players' learning dynamics are not completely independent since they track a common parameter together. Within the two-timescale learning framework, \cite{ref:Sayin20} presented independent learning dynamics for stochastic games and analyzed its almost-sure convergence also in two-player zero-sum stochastic games. Note that \cite{ref:OzdaglarICM} provides an overview of studies on independent learning dynamics in stochastic games.

Heterogenous learning in games has also been studied, however, with a specific focus on multi-timescale learning, e.g., \citep{ref:Leslie03,ref:Leslie05} for strategic-form games with repeated play and recently \citep{ref:Daskalakis20} for two-player zero-sum stochastic games. Particularly, dynamics of different players evolve at different timescales, e.g., one player's dynamics evolve slower than the others, contrary to \citep{ref:Leslie20,ref:Sayin20} where the two-timescale framework is at player-level. On heterogenous rates that may not lead to different timescales, \cite{ref:Zhu2011} studies heterogenous learning in a special case of zero-sum stochastic games while \cite{ref:Chasnov19} studies heterogenous gradient-based learning dynamics in continuous games.

In this paper, we address the heterogeneity and independence of learning in stochastic games by characterizing the convergence properties of the independent learning dynamics presented in \citep{ref:Sayin20} with player-dependent learning rates. This dynamics is a new variant of fictitious play combining the classical fictitious play \citep{ref:Fudenberg98} with the $Q$-learning \citep{ref:Watkins92} while they evolve at two different timescales. The key idea is that the underlying stochastic game can be viewed as a collection of auxiliary stage-games specific to each state whose payoff functions are the $Q$-functions. Though these auxiliary stage-games are not necessarily stationary, the slow evolution of $Q$-function estimates make them relatively stationary. The key challenge is the deviation of these auxiliary stage-games from the zero-sum structure due to the independent update of the $Q$-function estimates, and heterogenous learning rates boost this deviation further. We show the almost sure convergence of the dynamics under the usual two-timescale stochastic approximation assumptions when the discount factor is less than the product of the ratios of player-dependent step sizes. We elaborate on the implications and high-level interpretation of this result later in Section \ref{sec:result}. To show this result, we formulate a novel Lyapunov function formulation which reduces to the one presented in \citep{ref:Sayin20} in the homogenous case, and present a one-sided asynchronous convergence result, which has a similar flavor with \cite[Theorem 1]{ref:Tsitsiklis94}. 

We emphasize that the dynamics presented is different from equilibrium computation methods such as Shapley's value iteration \citep{ref:Shapley53} or its model-free version Minimax-Q algorithm \citep{ref:Littman94} by not requiring players to know the opponent's objective, i.e., the zero-sum structure of their stage-payoffs. Our dynamics also differs from multi-timescale learning schemes such as \citep{ref:Leslie03,ref:Daskalakis20} by addressing possible heterogeneity in (comparable) learning rates of players rather than exploiting it to characterize its convergence behavior especially since a convergence result with common learning rates is challenging. Furthermore, \cite{ref:Daskalakis20} considers the repeated play of a stochastic game with indefinite termination time allowing players to revise and improve their policies across repetitions. On the other hand, here we focus on the convergence of beliefs within a single stochastic game played over infinite horizon.

The rest of the paper is organized as follows. In Section \ref{sec:model}, we introduce stochastic games and describe the learning dynamics. We present the assumptions and main convergence result in Section \ref{sec:result} and the proof of the main convergence result in Section \ref{sec:proof}. In Section \ref{sec:example}, we provide an illustrative example. We conclude the paper with some remarks in Section \ref{sec:conclusion}. Appendices \ref{app:DI}-\ref{app:bounds} include the proofs of the technical lemmas used in Section \ref{sec:proof}.


\section{Independent Learning in (Zero-sum) Stochastic Games}\label{sec:model}

Formally, a two-player stochastic game is characterized by a tuple $\langle S,A,r^1,r^2,p,\gamma\rangle$.\footnote{For easy referral, we index players as player $1$ and player $2$. Furthermore, player $i$ is the typical player and player $-i$ is her opponent.} The \textit{finite} set of states is denoted by $S$ while $A=A^1\times A^2$ with $A^i$ denoting the \textit{finite} set of actions that player $i$ can take at any state.\footnote{The formulation can be extended to state-variant action sets straightforwardly.} The \textit{stage payoff function} of player $i$ is denoted by $r^i:S\times A \rightarrow \mathbb{R}$. In zero-sum case, we have $r^1(s,a) + r^2(s,a) = 0$ for all $(s,a)\in S\times A$. At any stage $k=0,1,\ldots$, if players play the action profile $a=(a^1,a^2)\in A$, then the state of the game, $s$, transits to another state, $s'$, according to the transition probability $p(s'|s,a)$. Player $i$'s objective is to maximize her expected sum of discounted stage-payoffs with the discount factor $\gamma \in [0,1)$.\footnote{Without a common discount factor, the underlying stochastic game may not be a zero-sum game even when $r^1(s,a) + r^2(s,a) = 0$ for all $(s,a)\in S\times A$.}

\cite{ref:Shapley53} (and \cite{ref:Fink64}) showed that in two-player zero-sum (and $n$-player general-sum) stochastic games, there always exists a \textit{stationary} equilibrium where players play \textit{stationary} strategies depending only on the current state. Let $\pi^i:S\rightarrow\Delta(A^i)$ denote the stationary strategy of player $i$ such that $\pi^i(s) \in \Delta(A^i)$ corresponds to her mixed strategy at state $s$, and $\pi = (\pi^1,\pi^2)$ denote the (stationary) strategy profile of players.\footnote{Given a set $A$, we denote the probability simplex over $A$ by $\Delta(A)$.} Then, the expected discounted sum of stage payoff of player $i$ under the strategy profile $\pi$ is given by
\be\label{eq:utility}
U^i(\pi) := \mathbb{E}\left\{\sum_{k=0}^{\infty} \gamma^k r^i(s_k,a_k)\right\},
\ee
where $a_k\sim \pi(s_k)$ denotes the action profile at stage $k$ while $\{s_k\}_{k\geq 0}$ is a stochastic process such that $s_k$ represents the state at stage $k$. The expectation is taken with respect to all randomness.

\begin{definition}[Stationary Nash Equilibrium]
We say that a stationary strategy profile $\pi_*$ is a stationary mixed-strategy equilibrium of the two-player stochastic game provided that
\be\label{eq:1}
U^i(\pi_*)\geq U^i(\pi^{i},\pi^{-i}_*), \quad \forall \pi^i\mbox{ and }i=1,2.
\ee
\end{definition}

We consider the same (independent) learning dynamics presented in \citep{ref:Sayin20} but with player-dependent learning rates to examine the robustness of the convergence result to such heterogeneity. Particularly, we can view the stage-wise interaction among players as they are playing \textit{auxiliary stage-games} specific to current state. For example, if player $i$ knew that the opponent will play according to the stationary strategy $\pi^{-i}$ in future stages, then her payoff function in the auxiliary stage-game specific to state $s$, denoted by $Q^i(s,\cdot):A\rightarrow\mathbb{R}$, would satisfy the following fixed-point condition
\be\label{eq:Q}
Q^i(s,a) = r^i(s,a) + \gamma \sum_{s'\in S} p(s'|s,a) \max_{\tilde{a}^i\in A^i} \mathbb{E}_{\tilde{a}^{-i}\sim\pi^{-i}(s')} \{Q^i(s',\tilde{a})\},\quad\forall a\in A.
\ee
This follows from backward induction that player $i$ would always look for maximizing her utility, as described in \eqref{eq:utility}. Note that the dependence on $\pi^{-i}$ is implicit for notational convenience. The function $Q^i(\cdot)$ is known as \textit{$Q$-function} in the MDP or reinforcement learning literature \citep{ref:Filar97book,ref:Sutton18}, and it is well-defined due to the contraction property of the Bellman operator. Therefore, the auxiliary stage game is the tuple $\langle A, Q^1(s,\cdot),Q^2(s,\cdot)\rangle$.

Players form a belief about their $Q$-function and adopt fictitious play in auxiliary stage-games by forming also a belief on the opponent strategy. We denote the beliefs of player $i$ at stage $k$ about the opponent strategy by $\hat{\pi}_k^{-i}$ and about her $Q$-function by $\hat{Q}_k^i$. At each stage $k$, player $i$ always take the best response action in the auxiliary stage-game based on her beliefs $\hat{\pi}_k^{-i}(s_k)$ and $\hat{Q}_k^i(s_k):=\hat{Q}_k^i(s_k,\cdot)$. Therefore, her action $a_k^i\in A^i$ always satisfies
\begin{equation}\label{eq:bestresponse}
a_k^i \in \argmax_{a^i\in A^i} \mathbb{E}_{a^{-i}\sim \hat{\pi}_k^{-i}(s_k)}\{\hat{Q}_k^i(s_k,a)\}.
\end{equation}
Player $i$ can observe the opponent's action $a_k^{-i}\in A^{-i}$ and update her beliefs according to 
\begin{subequations}\label{eq:update}
\begin{align}
&\hat{\pi}_{k+1}^{-i}(s) = \hat{\pi}_k^{-i}(s) + \mathbb{I}_{\{s=s_k\}}\alpha_{c_k(s)}^i (a_k^{-i} - \hat{\pi}_k^{-i}(s)),\label{eq:piupdate}\\
&\hat{Q}_{k+1}^{i}(s,a) = \hat{Q}_{k}^{i}(s,a) + \mathbb{I}_{\{s=s_k\}}\beta_{c_k(s)}^i\Big(r^i(s,a) + \gamma \sum_{s'\in S}p(s'|s,a) \hat{v}_k^i(s') - \hat{Q}_{k}^i(s,a)\Big),\label{eq:Qupdate}
\end{align}
\end{subequations}
for all $(s,a)$, where $\mathbb{I}_{\{s=s_k\}}$ is the indicator function, we let the pure action $a_k^{-i}$ be a deterministic strategy in the probability simplex $\Delta(A^{-i})$, the value function estimate $\hat{v}_k^i:S\rightarrow\mathbb{R}$ is defined by
\be\label{eq:v}
\hat{v}_k^i(s) = \max_{a^i\in A^i} \mathbb{E}_{a^{-i}\sim \hat{\pi}_k^{-i}(s)}\left\{\hat{Q}_k^i(s,a)\right\},
\ee 
and the player-dependent step sizes $\alpha_c^i\in(0,1)$ and $\beta_c^i\in (0,1)$ vanish with $c_k(s)$, the number of times state $s$ is visited until stage $k$. 

\section{Convergence Results}\label{sec:result}

The dynamics in \eqref{eq:update} allow player-dependent and belief-dependent step sizes.
In this section, we identify the conditions under which such a heterogenous two-timescale learning dynamics is guaranteed to converge to an equilibrium of the underlying two-player zero-sum stochastic game.

\begin{assumption}[Markov Chain]\label{assume:Markov}
Every state is visited infinitely often with probability $1$.
\end{assumption}

Beliefs associated with a state gets updated only if that state is visited. This assumption ensures that beliefs associated with each state gets updated infinitely often. Importantly, state transitions are controlled by players while players take actions according to their beliefs. Without Assumption \ref{assume:Markov}, it might be possible that they form incorrect beliefs about the value of other states such that these states would not get visited due to their greedy best responses (though they could have been visited for some non-greedy action). Then, the players are not able to revise and improve their beliefs about the values of these states. Therefore, their beliefs may converge to a \textit{self-confirming equilibrium} (a concept introduced for learning in extensive-form games to address out-of-equilibrium actions \citep{ref:Fudenberg95}) rather than a Nash equilibrium. We leave it as a future research direction.

Furthermore, Assumption \ref{assume:Markov} holds if the underlying stochastic game is \textit{irreducible}, e.g., transition probabilities between any pair of states are positive for any joint action as in \cite{ref:Leslie20}. This can be a restrictive assumption in practical application when players take deterministic best responses. However, it can be relaxed further as discussed in \cite{ref:OzdaglarICM} if players choose strategies in which every action is taken with some positive probability, e.g., due to smoothed best response or exploration, especially in the model-free cases where players do not know the stage-payoff and state transition kernel. We leave them as future research directions.

\begin{assumption}[Step Sizes]\label{assume:step}
The step sizes satisfy the following conditions:
\begin{itemize}
\item[$(a)$] For each $i=1,2$, the step sizes $\alpha_c^i\rightarrow 0$, $\beta_c^i\rightarrow 0$, and $\beta_c^i/\alpha_c^i\rightarrow0$ as $c\rightarrow\infty$, and the series $\sum_{c=0}^{\infty}\alpha_c^i = \sum_{c=0}^{\infty}\beta_c^i =\infty$.
\item[$(b)$]  Let $i,j\in\{1,2\}$ denote arbitrary player indices with $i\neq -i$ and $j\neq -j$. Then, for some $d_{\alpha},d_{\beta}\in (0,1]$, we have
$\lim_{c\rightarrow \infty}\alpha_c^i/\alpha_c^{-i} = d_{\alpha}$ and $\lim_{c\rightarrow \infty}\beta_c^j/\beta_c^{-j} = d_{\beta}$.\footnote{This is without loss of generality because if $\lim_{c\rightarrow\infty} \alpha_c^1/\alpha_c^2 \geq 1$ then we have $\lim_{c\rightarrow\infty}\alpha_c^2/\alpha_c^1 \leq 1$.}
\end{itemize} 
\end{assumption}

In Assumption \ref{assume:step}, part-$(a)$ is standard in multi-timescale stochastic approximation, e.g., see \citep[Chapter 6]{ref:Borkar08book}. 
Note that we do not need the assumption that the step sizes are square-summable because there is no stochastic approximation error (e.g., that can be induced from sampling noise) in the differential inclusion approximation of the discrete-time update rule. The part-$(b)$ says that the ratios of the step sizes have a non-zero limit at each timescale, i.e., players use comparable learning rates with each other. The limit assumption can be relaxed into conditions on limit inferior and limit superior, necessitating more involved analysis. We leave it as a future research direction.

\begin{theorem}[Convergence Result]\label{thm:main}
Given a two-player zero-sum stochastic game, suppose that players follow the heterogenous and independent learning dynamics \eqref{eq:update}. Under Assumptions \ref{assume:Markov} and \ref{assume:step}, we have 
$\hat{\pi}_{k}^i\rightarrow \pi_*^i$ and $\hat{Q}_k^i\rightarrow Q_*^i$ for each $i=1,2$ as $k\rightarrow \infty$ with probability $1$, for some stationary equilibrium $\pi_*$ and the associated $Q$-functions $(Q_*^1,Q_*^2)$ provided that
$\gamma \leq d_{\alpha}d_{\beta}$.
\end{theorem}

We attribute the upper bound on the discount factor as the heterogeneity increases how much the auxiliary stage-games can deviate from the zero-sum structure. However, a small discount factor will compensate this by restraining the deviation since the stage-payoffs have zero-sum structure. Furthermore, when players have (asymptotically) common step sizes, i.e., $d_{\alpha}=d_{\beta}=1$, the bound on the discount factor becomes $\gamma < 1$, which is inherent to the discounted stochastic games. Therefore, Theorem \ref{thm:main} reduces to the convergence result \cite[Theorem 4.3]{ref:Sayin20} for the special case of homogeneous learning rates. 

We can interpret Theorem \ref{thm:main} at a high level as players could reach to an equilibrium through heterogenous and independent learning dynamics if they are sufficiently myopic so that they discount the impact of future stages more in their utilities. We can also interpret the discount factor as the continuation probability of the stochastic game with indefinite horizon length \citep{ref:Shapley53}. Therefore, the heterogenous learning dynamics is guaranteed to converge to an equilibrium if the game has sufficiently short expected termination time.

Note also that stochastic games turn into strategic-form games with repeated play if there is only one state and $\gamma=0$. In that case, the assumption on $\gamma$ is always satisfied and the auxiliary stage-games are always zero-sum. Therefore, we have the following corollary to Theorem \ref{thm:main}.

\begin{corollary}[Heterogenous Fictitious Play]
Consider a two-player zero-sum strategic-form game played repeatedly. Suppose that players follow the fictitious play dynamics with player-dependent learning rates $\alpha_k^1,\alpha_k^2$ such that their ratio has non-zero limit. Then, the beliefs formed about the opponent strategy converge to an equilibrium of the game.
\end{corollary}


\section{Proof of the Convergence Result - Theorem \ref{thm:main}}\label{sec:proof}

In two-player zero-sum stochastic games, \cite{ref:Shapley53} provided a (minimax) value iteration to compute equilibrium values associated with a stationary equilibrium. The operator used in Shapley's value iteration can be transformed into
\be\label{eq:operator}
(\mathcal{F}^iQ^i)(s,a) = r^i(s,a) + \gamma \sum_{s'\in S} p(s'|s,a) \val^i(Q^i(s',\cdot)),\quad\forall (s,a)\in S\times A,
\ee
as in \citep{ref:Szepesvari99}, where the minimax value function $\val^i(\cdot)$ is defined by
\be
\val^i(Q^i(s,\cdot)) := \max_{\mu^i\in\Delta(A^i)}\min_{\mu^{-i}\in\Delta(A^{-i})}\mathbb{E}_{(a^i,a^{-i})\sim (\mu^i,\mu^{-i})} \{Q^i(s,a)\}.
\ee
The operator $\mathcal{F}^i$ is a contraction and its unique fixed point $Q_*^i = \mathcal{F}^iQ_*^i$ is the equilibrium $Q$-function of the underlying stochastic game. 

In this proof, we look for identifying the conditions under which 
\be\label{eq:e}
e_k^i(s) := \hat{v}_k^i(s) - \val^i(\hat{Q}^i_k(s))\rightarrow 0\quad\mbox{and}\quad \tilde{Q}_k(s,a):= \hat{Q}_k^i(s,a)-Q_*^i(s,a)\rightarrow 0,
\ee
for each $(s,a)$ and $i=1,2$, by using stochastic differential inclusion theory while formulating a novel Lyapunov function and one-sided convergence result to address the heterogeneity. To this end, firstly, the following lemma establishes the connection between \eqref{eq:update} and its (continuous-time) limiting differential inclusion based on \citep{ref:Benaim05}.\footnote{Without loss of generality, suppose that $\lim_{c\rightarrow\infty} \alpha_c^1/\alpha_c^2 = d_{\alpha}\in (0,1]$.} 
The proof is deferred to Appendix \ref{app:DI}. 

\begin{lemma}\label{lem:DI}
For each state $s$, the limiting differential inclusion of \eqref{eq:update} is given by
\begin{align}\label{eq:DI}
&\dot{\pi}^1 + d_{\alpha} \pi^1 \in d_{\alpha} \argmax_{a^1\in A^1} \mathbb{E}_{a^2\sim\pi^2}\{Q^1(a)\}\quad\mbox{and}\quad\dot{\pi}^2 + \pi^2 \in \argmax_{a^2\in A^2} \mathbb{E}_{a^1\sim\pi^1}\{Q^2(a)\}
\end{align}
and $\dot{Q}^1(a) = 0$ for all $(i,a)$, where $\pi^i:[0,\infty)\rightarrow \Delta(A^i)$ and $Q^i(a):[0,\infty)\rightarrow \mathbb{R}$, for each $i=1,2$, are continuous-time functions where we drop the dependence on $s$ for notational convenience.
\end{lemma}

We note that there always exists an absolutely continuous solution to \eqref{eq:DI} since the best response satisfies the conditions listed in \cite[Hypothesis 1.1]{ref:Benaim05}. Then, \cite[Theorem 3.6 and Proposition 3.27]{ref:Benaim05} yield that we can characterize the convergence properties of the discrete-time dynamics \eqref{eq:update} in terms of the zero-set of a Lyapunov function to the differential inclusion \eqref{eq:DI}. 

Though \eqref{eq:DI} resembles to continuous-time best response dynamics in the auxiliary stage-game with time-invariant payoff functions $Q^1$ and $Q^2$, there are two challenges: $(i)$ the heterogeneity when $d_{\alpha}\in(0,1)$ and $(ii)$ the deviation from the zero-sum structure since $Q_k^1(s,a)+Q_k^2(s,a)$ is not necessarily zero for all $(s,a)$ and $k$ when they are updated independently according to \eqref{eq:Qupdate}.  Our candidate Lyapunov function is defined by
\be\label{eq:Lyapunov}
V(\pi,Q) := \left(d_{\alpha}\Delta^1(\pi^2,Q^1) + \Delta^2(\pi^1,Q^2) - \Xi(Q^1,Q^2)\right)_+,
\ee
where we define 
\begin{align}
&\Delta^i(\pi^{-i},Q^i):= \max_{a^i\in A^i}\mathbb{E}_{a^{-i}\sim\pi^{-i}}\{Q^i(a)\} - \val^i(Q^i) \geq 0\label{eq:Delta}\\
&\Xi(Q^1,Q^2) := \lambda\|Q^1+Q^2\| - (\val^1(Q^1) + \val^2(Q^2)),
\end{align}
where $\lambda\in(1,d_{\alpha}d_{\beta}/\gamma)$,  $\|\cdot\|$ is the maximum norm, i.e., $\|Q\| = \max_{a}|Q(a)|$ and the positive function $(x)_+=\max\{0,x\}$. The candidate function is non-negative by its definition. Its zero-set $\{(\pi,Q):V(\pi,Q)=0\}$ is given by
\be\label{eq:zero}
\left\{(\pi,Q): \Delta^1(\pi^2,Q^1) + \Delta^2(\pi^1,Q^2) \leq \Xi(Q^1,Q^2) + (1-d_{\alpha})\Delta^1(\pi^2,Q^1) \right\}
\ee
and $\Delta^i$ is the continuous-time counterpart of the tracking error \eqref{eq:e}. Therefore, the convergence to the zero-set \eqref{eq:zero} would provide an (asymptotic) upper bound on the sum of tracking errors.
Note also that when $d_{\alpha}=1$, minimax values disappear and $V(\cdot)$ reduces to the one presented in \citep{ref:Sayin20}. 

Since $Q^i$'s are time-invariant, $\val^i(Q^i)$'s and $\Xi$ are also time-invariant. As shown in the following lemma, these time-invariant terms together with the positive function play an important role for the validity of the candidate $V(\cdot)$ as a Lyapunov function to \eqref{eq:DI} for the zero-set \eqref{eq:zero} and later in characterizing the convergence properties of the tracking error. 
The proof is deferred to Appendix \ref{app:Lyapunov}.\footnote{The condition that $r^1(s,a)+r^2(s,a) = 0$ for all $(s,a)$ plays an important role in the proof of Lemma \ref{lem:Lyapunov}.}

\begin{lemma}[Lyapunov Function]\label{lem:Lyapunov}
The candidate function $V(\cdot)$ is a Lyapunov function of the differential inclusion \eqref{eq:DI} for the zero-set $\{(\pi,Q):V(\pi,Q)=0\}$. In other words, for any absolutely continuous solution to \eqref{eq:DI}, we have $V(\pi(t'),Q(t'))<V(\pi(t),Q(t))$ for all $t'>t$ if $V(\pi(t),Q(t)) > 0$, and $V(\pi(t'),Q(t')) = 0$ for all $t'>t$ if $V(\pi(t),Q(t)) = 0$.
\end{lemma}

Based on the stochastic differential inclusion theory  \cite{ref:Benaim05}, Lemma \ref{lem:Lyapunov}, the definition of $\hat{v}_k^i$, as described in \eqref{eq:v}, and the zero-set \eqref{eq:zero}, we can conclude that
\be\label{eq:up}
(\bar{v}_k(s) -\lambda \|\bar{Q}_k(s)\| - (1-d_{\alpha})e_k^1(s))_+ \rightarrow 0,
\ee
where we define $\bar{v}_k(s):=\hat{v}_k^1(s)+\hat{v}_k^2(s)$ and $\bar{Q}_k(s) := \hat{Q}_k^1(s)+\hat{Q}_k^2(s)$ with the tracking error $e_k^1$, as described in \eqref{eq:e}. Since we let $i=1$ arbitrarily for notational convenience, \eqref{eq:up} can be written as
\be\label{eq:upp}
\bar{v}_k(s) \leq \lambda \|\bar{Q}_k(s)\| + (1-d_{\alpha})e_k^i(s) + \epsilon_k(s),
\ee
where $\epsilon_k(s)\rightarrow 0$ as $k\rightarrow\infty$ almost surely. The upper bound on $\bar{v}_k(s)$ is in terms of $\|\bar{Q}_k(s)\|$ and the tracking error $e_k^i$. The following lemma provides an upper bound on the tracking error so that we can obtain an upper bound in terms of $\bar{Q}_k(s)$ only. The proof is deferred to Appendix \ref{app:aux}.

\begin{lemma}\label{lem:aux}
We have $0\leq e_k^i(s) = \hat{v}_k^i(s)-\val^i(\hat{Q}_k^i(s)) \leq \bar{v}_k(s) - \min_{a\in A}\,\bar{Q}_k(s,a)$ for every $s$.
\end{lemma}

Based on \eqref{eq:upp} and Lemma \ref{lem:aux}, we obtain
\begin{align}\label{eq:belabo}
 \min_{a\in A} \, \bar{Q}_k(s,a) \leq \bar{v}_k(s) \leq \frac{\lambda}{d_{\alpha}} \|\bar{Q}_k(s)\| - \frac{1-d_{\alpha}}{d_{\alpha}} \min_{a\in A}\,\bar{Q}_k(s,a) + \frac{1}{d_{\alpha}} \epsilon_k(s),
\end{align}
where the lower bound follows since $\hat{v}_k^i(s) \geq \mathbb{E}_{a\sim\hat{\pi}_k(s)}\{\hat{Q}_k^i(s,a)\}$ for $i=1,2$. We emphasize that in the homogenous case,  the evolution of $\bar{Q}_k(s,a)$ is given by
\be
\bar{Q}_{k+1}(s,a) = \bar{Q}_k(s,a)+\mathbb{I}_{\{s=s_k\}}\beta_{c_k(s)} \left(\gamma \sum_{s'\in S}p(s'|s,a) \bar{v}_k(s') - \bar{Q}_k(s,a)\right)
\ee 
due to the symmetry that $\beta_c^i = \beta_c$. However, the characterization of the convergence properties of $\bar{Q}_k(s,a)$ in the heterogenous case requires more involved analysis where we will address the limit inferior and limit superior of $\bar{Q}_k(s,a)$ separately.

Recall the definition of $\tilde{Q}_k^i$ in \eqref{eq:e} and note that the fixed points $Q_*^1$ and $Q_*^2$ satisfy $Q_*^1(s,a) + Q_*^2(s,a) = 0$ for every $(s,a)$. Therefore, we also have
$\bar{Q}_k= \tilde{Q}_k^1 + \tilde{Q}_k^2$. Based on \eqref{eq:Qupdate} and the definition of the fixed point $Q_*^i$, the evolution of $\tilde{Q}_k^i$ can be written as
\begin{align}\label{eq:tilde}
\tilde{Q}_{k+1}^i(s,a) = &\tilde{Q}_k^i(s,a) + \bar{\beta}_k^i(s)\left(\gamma \sum_{s'\in S}p(s'|s,a) (\hat{v}_k^i(s') - \val^i(Q_*^i(s'))) - \tilde{Q}_k^i(s,a)\right),
\end{align}
where $\bar{\beta}_k^i(s):=\mathbb{I}_{\{s=s_k\}}\beta_{c_k(s)}^i$ and $r^i(s,a)$ disappears since $Q_*^i(s,a)=(\mathcal{F}^iQ_*^i)(s,a)$. We are interested in the limit inferior of $\tilde{Q}_k^i$, for each $i$, so that we can formulate the limit inferior of $\bar{Q}_k$, which will play an important role in \eqref{eq:belabo}. 

The following lemma characterizes the limit inferior of an iterate whose evolution satisfies one-sided contraction-like condition. The proof is deferred to Appendix \ref{app:async}.

\begin{lemma}[One-sided Asynchronous Discrete-time Convergence]\label{lem:async}
Consider a sequence of vectors $\{y_k\}_{k=0}^{\infty}$ such that the $n$th entry, denoted by $y_k(n)$, satisfies the following lower bound:
\be
y_{k+1}(n) \geq y_k(n) + \beta_{k}(n)\left(\gamma \min_{m}y_k(m)-y_k(n) + \epsilon_k(n)\right),
\ee
where $\gamma \in (0,1)$, the vanishing (possibly random) step size $\beta_k(n)\in [0,1]$ satisfies $\beta_k(n)\rightarrow 0$ and $\sum_{k=0}^{\infty}\beta_k(n)=\infty$, and $\liminf_k \epsilon_k(n)\geq 0$  for each $n$, with probability $1$. Suppose that $\min_n y_k(n) \geq M$ for all $k$. Then, we have 
$\liminf_{k} y_k(n) \geq 0$ for all $n$, with probability $1$.
\end{lemma}

By the definition of the tracking error \eqref{eq:e}, we have 
\begin{align}
\sum_{s'\in S}p(s'|s,a) (\hat{v}_k^i(s') - \val^i(Q_*^i(s'))) &=  \sum_{s'\in S} p(s'|s,a) (e_k^i(s') + \val^i(\hat{Q}_k^i(s')) - \val^i(Q_*^i(s')))\nn
\end{align}
which is greater than $\min_{(s,a)} \tilde{Q}_k^i(s,a)$
since $e_k^i(s')\geq 0$ and
$\val^i(\hat{Q}_k^i(s')) - \val^i(Q_*^i(s'))\geq \min_{a\in A}\;\tilde{Q}_k^i(s',a)$
for all $s'$. Under Assumptions \ref{assume:Markov} and \ref{assume:step}, the step size $\bar{\beta}_k(s)\in[0,1]$ vanishes and $\sum_{k=0}^{\infty}\bar{\beta}_k(n) = \infty$ with probability. Therefore, we can invoke Lemma \ref{lem:async}, for every $i=1,2$, and obtain
\be
\liminf_{k\rightarrow\infty} \tilde{Q}_k^i(s,a) \geq 0\quad\Rightarrow\quad \liminf_{k\rightarrow\infty} \bar{Q}_k(s,a) \geq 0,\label{eq:liminf}
\ee
almost surely for every $(s,a)$. Combined with \eqref{eq:belabo}, the bound \eqref{eq:liminf} yields that
\be\label{eq:final}
\underline{\epsilon}_k(s) \leq \bar{v}_k(s) \leq \frac{\lambda}{d_{\alpha}} \|\bar{Q}_k(s)\| + \overline{\epsilon}_k(s),\quad\forall s,
\ee
for some error terms $\underline{\epsilon}_k(s)\rightarrow 0$ and $\overline{\epsilon}_k(s)\rightarrow 0$ as $k\rightarrow\infty$ for each $s$ almost surely. 

The heterogeneity of $\beta_c^i$'s has not played a role up to this point. Suppose that $\lim_{c\rightarrow\infty}\beta_c^1/\beta_c^2 = d_{\beta}\in (0,1]$ without loss of generality. Then, the iteration \eqref{eq:tilde} can be written as
\be\label{eq:entry}
\begin{bmatrix} \tilde{Q}_{k+1}^1\\ \tilde{Q}_{k+1}^2 \end{bmatrix} = 
\begin{bmatrix} \tilde{Q}_k^1 \\ \tilde{Q}_k^2\end{bmatrix} + \bar{\beta}_k(s) \begin{bmatrix} d_{\beta} \left(\gamma\sum_{s'}p(s'|\cdot)(\hat{v}_k^1(s') - \val^1(Q_*^1(s'))) - \tilde{Q}_k^1 + \zeta_k\right)\\ \gamma\sum_{s'}p(s'|\cdot)(\hat{v}_k^2(s') - \val^2(Q_*^2(s'))) - \tilde{Q}_k^2\end{bmatrix},
\ee
where we dropped the argument $(s,a)$ for notational convenience and the error term is given by
\be
\zeta_k(s,a) := \frac{1}{d_{\beta}}\left(\frac{\beta_{c_k(s)}^1}{\beta_{c_k(s)}^2}-d_{\beta}\right)\left(\gamma\sum_{s'}p(s'|s,a)(\hat{v}_k^1(s') - \val^1(Q_*^1(s'))) - \tilde{Q}_k^1(s,a)\right),
\ee
which is asymptotically negligible by Assumptions \ref{assume:Markov} and \ref{assume:step}, and by the boundedness of the iterates. For notational convenience, we also define
$\Gamma_k(s,a) := \tilde{Q}_k^1(s,a) + d_{\beta}\tilde{Q}_k^2(s,a)$.
Then, the weighted combination of the entries in \eqref{eq:entry} yield that
\be\label{eq:Gamma}
\Gamma_{k+1}(s,a) = \Gamma_k(s,a) + \bar{\beta}_k(s) d_{\beta}\left(\gamma\sum_{s'\in S}p(s'|s,a)\bar{v}_k(s') -\bar{Q}_k(s,a)\right)
\ee
since $\val^1(Q_*^1(s))+\val^2(Q_*^2(s))=0$ for each $s$. By \eqref{eq:liminf} and \eqref{eq:Gamma}, we also have
\be\label{eq:liminf2}
\liminf_{k\rightarrow\infty}\Gamma_k(s,a) \geq 0,\quad\forall\;(s,a)
\ee
with probability $1$. In order to characterize the convergence properties of $\Gamma_k$ within the framework of Lemma \ref{lem:async}, we introduce the following lemma formulating a bound on $\bar{Q}_k$ in terms of $\Gamma_k$. The proof is deferred to Appendix \ref{app:bounds}. 

\begin{lemma}\label{lem:bounds}
For some $\underline{\eta}_k(s,a)\rightarrow 0$ and $\overline{\eta}_k(s,a)\rightarrow0$ with probability $1$, we have
\be
\overline{\eta}_k(s,a) + \frac{1}{d_{\beta}}\Gamma_k(s,a) \geq |\bar{Q}_k(s,a)| \geq \bar{Q}_k(s,a) \geq \Gamma_k(s,a) + \underline{\eta}_k(s,a),\quad\forall (s,a).
\ee
\end{lemma}

Based on \eqref{eq:final}, \eqref{eq:Gamma}, and Lemma \ref{lem:bounds}, we obtain
\be
\Gamma_{k+1}(s,a) \leq \Gamma_k(s,a) + \bar{\beta}_k(s) d_{\beta}\left(\frac{\gamma\lambda}{d_{\alpha}d_{\beta}}\max_{(s',a')}\; \Gamma_k(s',a') - \Gamma_k(s,a) +\varepsilon_k(s,a)\right),
\ee
where $\varepsilon_k(s,a)$ is an asymptotically negligible error almost surely. We can invoke Lemma \ref{lem:async} for $\{-\Gamma_k\}_{k\geq 0}$ and obtain that
\be\label{eq:limsup}
\liminf_{k\rightarrow\infty} -\Gamma_k(s,a) = -\limsup_{k\rightarrow\infty} \Gamma_{k}(s,a) \geq 0\quad\Rightarrow\quad \limsup_{k\rightarrow\infty} \Gamma_{k}(s,a) \leq 0.
\ee
Combined with \eqref{eq:liminf2}, the bound \eqref{eq:limsup} yields that
$\lim_{k\rightarrow\infty}\Gamma_k(s,a) = 0$. Hence, $\lim_{k\rightarrow\infty}\tilde{Q}_k^i(s,a) = 0$, for all $(s,a)$ and $i=1,2$ with probability $1$. Since $\bar{Q}_k = \tilde{Q}_k^1+\tilde{Q}_k^2$, we also have $\bar{Q}_k(s,a)\rightarrow0$ for each $(s,a)$. Then, we obtain $\bar{v}_k(s)\rightarrow 0$ for each $s$ by \eqref{eq:final}. Therefore, Lemma \ref{lem:aux} yields that the tracking error $e_k^i(s)\rightarrow 0 $ for each $s$, with probability $1$.  


\begin{figure}[t!]
\centering
\includegraphics[width=.55\textwidth]{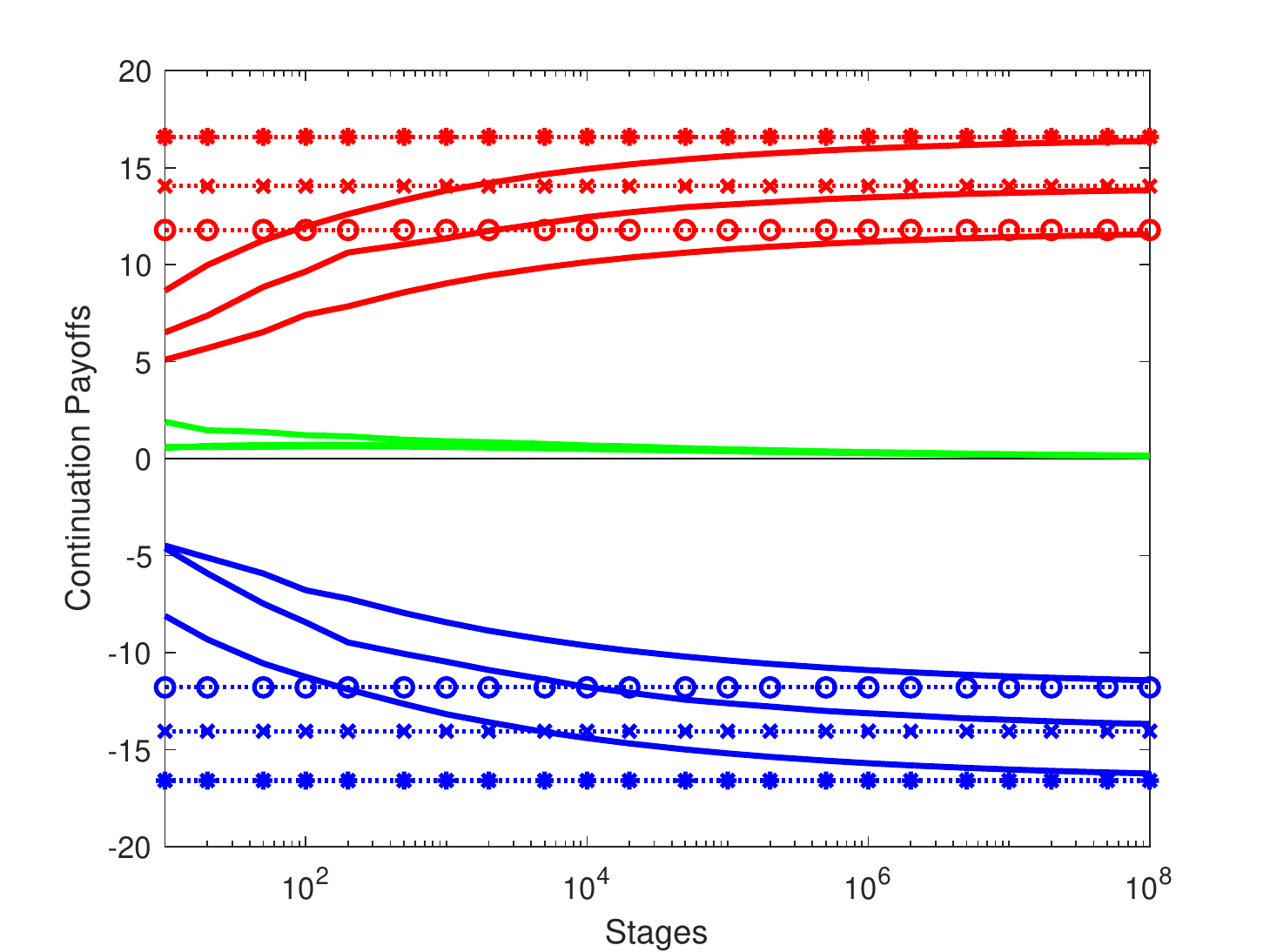}
\caption{Evolution of $\hat{v}_k^1$, $\hat{v}_k^2$, and $\bar{v}_k=v^1_k + v_k^2$ for each $s$, respectively, in {\color{red!80!black}\bf red}, {\color{blue!80!black}\bf blue}, and {\color{green!80!black} \bf green}. 
The dotted lines ${\color{red}\ldots}/{\color{blue}\ldots}$ are the actual equilibrium values associated with each state.} 
\label{fig:plot}
\end{figure}

\section{Numerical Example}\label{sec:example}

In this section, we examine the convergence properties of the heterogenous learning dynamics numerically in an (irreducible) zero-sum stochastic game whose configuration is selected arbitrarily as in \citep{ref:Sayin20}. For example, there are three states, four actions per state, and the discount factor is $0.8$. The player-dependent step sizes are set as $\alpha_c^1 = (1+c)^{0.5}$ and $\alpha_c^2 = (1+0.81\cdot c)^{0.5}$ while $\beta_c^1=(1+c)^{-1}$ and $\beta_c^2=(1+0.95\cdot c)^{-1}$. In Figure \ref{fig:plot}, we plot the evolution of the value function estimates of both player in addition to $\bar{v}_k$ to illustrate that the auxiliary stage-games become zero-sum also in the heterogenous case. Since the underlying game is irreducible, Assumption \ref{assume:Markov} holds. For these step sizes, we have $d_{\alpha} = 0.9$ and $d_{\beta} = 0.95$. Assumption \ref{assume:step} also holds. Therefore, Theorem \ref{thm:main} says that the dynamics should converge to an equilibrium of the game and we have observed the convergence of the value function estimates to the equilibrium values of the game, as expected from Theorem \ref{thm:main}. 

\section{Conclusion}\label{sec:conclusion}

We showed the almost sure convergence of two-timescale fictitious play with heterogeneous learning rates in two-player zero-sum stochastic games under the standard assumptions in two-timescale stochastic approximation methods when the discount factor is less than the product of the ratios of the player-dependent step-sizes. Since strategic-form games played repeatedly is a special case of stochastic games, this result also implied the almost sure convergence of fictitious play with heterogenous rates in zero-sum strategic-form games. We attributed the bound on the discount factor to the deviation of the auxiliary stage-games from the zero-sum structure, which becomes multifold with the player-dependent rates while sufficiently small discount rates can compensate it. Future research directions include characterizing the convergence properties of heterogenous and independent learning dynamics in stochastic games other than zero-sum, and in model-free and minimal information cases, e.g., as in \citep{ref:Sayin21}. 

\acks{This is an extended version with Appendices \ref{app:DI}-\ref{app:bounds} for the proofs of technical lemmas used in the proof of Theorem \ref{thm:main} in Section \ref{sec:proof}.}

\appendix

\section{Proof of Lemma \ref{lem:DI}}\label{app:DI}

The dynamics specific to state $s$, i.e., \eqref{eq:piupdate} and \eqref{eq:Qupdate}, can be written as
\be\label{eq:pi2}
\begin{bmatrix}\hat{\pi}_{k+1}^1(s) \\ \hat{\pi}_{k+1}^{2}(s)\end{bmatrix} = \begin{bmatrix}\hat{\pi}_k^1(s)\\\hat{\pi}_k^{2}(s)\end{bmatrix} + \bar{\alpha}_k(s)\left(\begin{bmatrix} d_{\alpha}(a_k^1 - \hat{\pi}_k^1(s)) \\ a_k^{2} - \hat{\pi}_k^{2}(s)\end{bmatrix} + \begin{bmatrix} \varepsilon_{k}(s) \\ \mathbf{0} \end{bmatrix}\right),
\ee
and
\be\label{eq:Q2}
\begin{bmatrix} \hat{Q}_{k+1}^1(s,a) \\ \hat{Q}_{k+1}^{2}(s,a) \end{bmatrix} = \begin{bmatrix} \hat{Q}_{k}^1(s,a) \\ \hat{Q}_{k}^{2}(s,a) \end{bmatrix} + \bar{\alpha}_k(s)\begin{bmatrix} E_k^1(s,a) \\ E_k^2(s,a) \end{bmatrix},
\ee
where $\mathbf{0}$ is a zero vector, the step size $\bar{\alpha}_k(s) := \mathbb{I}_{\{s=s_k\}}\alpha_{c_k(s)}^{2}\in[0,1]$, and the error terms are defined by
\begin{subequations}\label{eq:error}
\begin{align}
&\varepsilon_{k}(s) := \left(\frac{\alpha_{c_k(s)}^1}{\alpha_{c_k(s)}^{2}} - d_{\alpha}\right)(a_k^1 - \hat{\pi}_k^1(s))\\
&E_k^1(s,a) := \frac{\alpha_{c_k(s)}^1}{\alpha_{c_k(s)}^2} \frac{\beta_{c_k(s)}^1}{\alpha_{c_k(s)}^1}\left(r^1(s,a) + \gamma \sum_{s'\in S} p(s'|s,a) \hat{v}_k^1(s') - \hat{Q}_k^1(s,a)\right) \\
&E_k^2(s,a) := \frac{\beta_{c_k(s)}^2}{\alpha_{c_k(s)}^2}\left(r^2(s,a) + \gamma \sum_{s'\in S} p(s'|s,a) \hat{v}_k^2(s') -\hat{Q}_k^2(s,a)\right).
\end{align}
\end{subequations}
Note that the iterates are bounded since $\gamma\in (0,1)$, the stage-payoffs have compact support and the step sizes are in $(0,1]$. Therefore, Assumptions \ref{assume:Markov} and \ref{assume:step} yield that the error terms \eqref{eq:error} are all asymptotically negligible for all $(s,a)$. Note also that under Assumption \ref{assume:Markov} and \ref{assume:step}, we have
$
\sum_{k=0}^{\infty}\bar{\alpha}_k(s) = \infty
$
while $\bar{\alpha}_k(s)\rightarrow 0$ as $k\rightarrow \infty$ with probability $1$. Furthermore, the best response satisfies the conditions for the stochastic differential inclusion theory \citep[Hypothesis 1.1]{ref:Benaim05}. Therefore, its limiting differential inclusion is given by \eqref{eq:DI}.

\section{Proof of Lemma \ref{lem:Lyapunov}}\label{app:Lyapunov}

For fixed absolutely continuous solution $(\pi^1(t),\pi^2,Q^1,Q^2)$ to \eqref{eq:DI}, the argument of the positive function in \eqref{eq:Lyapunov} is given by
\be
L(t) := d_{\alpha}\Delta^1(\pi^2(t),Q^1) + \Delta^2(\pi^1(t),Q^2) - \Xi(Q^1,Q^2),
\ee
which is also an absolutely continuous function since $\max$ and addition satisfy the Lipschitz condition \cite[Lemma 4.3.2]{ref:Bogachev20book}. Therefore, we can compute its derivative almost everywhere as in \cite{ref:Harris98} and obtain
\begin{align}
\frac{d}{dt}\left(d_{\alpha}\Delta^1+\Delta^2-\Xi\right) &= d_{\alpha}(a^1)^TQ^1\frac{d\pi^2}{dt} + \left(\frac{d\pi^1}{dt}\right)^TQ^2a^2\nn\\
&= d_{\alpha}\left((a^1)^TQ^1(a^2-\pi^2) + (a^1-\pi^1)^TQ^2a^2\right)\nn\\
&=d_{\alpha}\left((a^1)^T(Q^1+Q^2)a^2 - (\Delta^1+\Delta^2) - (\val^1(Q^1) + \val^2(Q^2))\right),\label{eq:derivative}
\end{align}
where
\be
a^1 := \frac{1}{d_{\alpha}}\frac{d\pi^1}{dt} + \pi^1\quad\mbox{and} \quad a^2 := \frac{d\pi^2}{dt}+\pi^2.
\ee
Note that if it is zero-sum, we have $\dot{L} = -d_{\alpha}(\Delta^1+\Delta^2)\leq 0$ and equal to zero only if $\Delta^i=0$, which corresponds to the equilibrium. Suppose that it is not zero-sum. Then, the time derivative of $L$ is not necessarily negative almost everywhere. However, we can write \eqref{eq:derivative} as
\begin{align}
\frac{d}{dt}\left(d_{\alpha}\Delta^1+\Delta^2-\Xi\right) =&\; -d_{\alpha}\overbrace{(d_{\alpha}\Delta^1+\Delta^2-\Xi)}^{=L} + d_{\alpha}((a^1)^T(Q^1+Q^2)a^2 - \lambda \|Q^1+Q^2\|)\nn\\
&- d_{\alpha}(1-d_{\alpha})\Delta^1.\label{eq:dev2}
\end{align}
The second term at the right-hand side is negative since $\lambda > 1$ and it is not zero-sum. Furthermore, the last term in \eqref{eq:dev2} is non-positive since $d_{\alpha}\in(0,1]$ and $\Delta^1\geq 0$. Therefore, we have 
\be
\dot{L} < -d_{\alpha} L, 
\ee
almost everywhere. Therefore, the absolutely continuous $L(t)$ is strictly decreasing when $L(t)\geq 0$ and $(-\infty,0]$ is a positively invariant set for $L(t)$. Correspondingly, $V(\pi(t),Q) = (L(t))_+$ is strictly decreasing when $V(\pi(t),Q)> 0$ and $(-\infty,0]$ is also a positively invariant set for $V(\pi(t),Q)$. 

Note that we can pick $\lambda$ such that $\lambda\in(1,d_{\alpha}d_{\beta}/\gamma)$ since it is arbitrary and $\gamma < d_{\alpha}d_{\beta}$, as stated in Theorem \ref{thm:main}. This completes the proof.

\section{Proof of Lemma \ref{lem:aux}}\label{app:aux}

The proof follows from the saddle point inequality:
\begin{align}
\hat{v}_k^i(s) = \max_{a^i\in A^i} \mathbb{E}_{a^{-i}\sim\hat{\pi}_k^{-i}(s)} \{\hat{Q}_k^i(s,a)\} \geq \val^i(\hat{Q}_k^i(s)) \geq \min_{a^{-i}\in A^{-i}} \mathbb{E}_{a^i\sim\hat{\pi}_k^i(s)}\{\hat{Q}_k^i(s)\}
\end{align}
because the right-most term is bounded from below by
\begin{align}
\min_{a^{-i}} \mathbb{E}_{a^i\sim\hat{\pi}_k^i(s)}\{\bar{Q}_k(s) - \hat{Q}_k^{-i}(s)\} &\geq \min_{a^{-i}} \mathbb{E}_{a^i\sim\hat{\pi}_k^i(s)}\{\bar{Q}_k(s)\} + \min_{a^{-i}} \mathbb{E}_{a^i\sim\hat{\pi}_k^i(s)}\{- \hat{Q}_k^{-i}(s)\}\nn\\
&=\min_{a^{-i}} \mathbb{E}_{a^i\sim\hat{\pi}_k^i(s)}\{\bar{Q}_k(s)\} - \max_{a^{-i}} \mathbb{E}_{a^i\sim\hat{\pi}_k^i(s)}\{\hat{Q}_k^{-i}(s)\}.
\end{align}
Therefore, we obtain
\be
\hat{v}_k^i(s) \geq \val^i(\hat{Q}_k^i(s)) \geq -\hat{v}_k^{-i}(s) + \min_{a\in A}\,\bar{Q}_k(s,a).
\ee
The difference between the first and the second term is bounded from above by the difference between the first and the third term. This completes the proof.

\section{Proof of Lemma \ref{lem:async}}\label{app:async}

Note that if $M\geq 0$, the lower bound, $y_k(n)\geq M$ for all $k\geq 0$ and $n$, already implies that $\liminf_{k\rightarrow\infty}y_k(n)\geq 0$ for all $n$. 

Suppose that $M<0$. Then, the proof has a flavor similar to \cite[Theorem 1]{ref:Tsitsiklis94} and \cite[Theorem 5.1]{ref:Sayin20} but only for one-side to characterize the limit inferior of the sequence based on the assumption that all iterates are bounded from below and update has a one-sided contraction-like structure. 

Define the negative sequence $\{M^t< 0\}_{t\geq0}$ over a separate timescale by
\be
M^{t+1} = (\gamma + 2\epsilon) M^t,\quad\forall\;t\geq 0,
\ee
and $M^0 = M$, where $\epsilon \in (0,(1-\gamma)/2)$.  Since $\gamma+2\epsilon \in(0,1)$, we have $M^t\rightarrow 0$ monotonically (from below) as $t\rightarrow\infty$.

Since $\liminf_k\epsilon_k(n)\geq0$ for each $n$ and $y_k(n)\geq M^0$ for all $k\geq 0$, there exists $k^0$ such that 
\be
y_k(n)\geq M^0,\quad\mbox{and}\quad\epsilon_k(n)>\epsilon M^0\quad\forall \;k\geq k^0.
\ee 
Suppose that for some $t\geq 0$, there exists $k^t$ such that 
\be
y_k(n)\geq M^t,\quad\mbox{and}\quad\epsilon_k(n)>\epsilon M^t\quad\forall \;k\geq k^t.
\ee 
Then, we can define an auxiliary sequence $\{Y_k^t\}_{k\geq k^t}$ by
\be
Y_{k+1}^t(n) = Y_k^t(n)(1-\beta_k(n)) + \beta_k(n) (\gamma + \epsilon)M^t,\quad\forall\;k\geq k^t,
\ee
and $Y_{k^t}^t(n) = M^t$, for each $n$ such that
\be
Y_k^t(n) \leq y_k(n),\quad\forall\; k\geq k^t, 
\ee
for each $n$, by its definition. Note that $Y_k^t(n)\rightarrow (\gamma+\epsilon) M^t$ for each $n$ as $k\rightarrow\infty$ almost surely. Since $(\gamma+2\epsilon)M^t<(\gamma+\epsilon)M^t$, there exists $k^{t+1}\geq k^t$ such that 
\be
y_k(n)\geq (\gamma+2\epsilon)M^t=M^{t+1}\quad\mbox{and}\quad \epsilon_k(n)>\epsilon M^{t+1},\quad\forall k\geq k^{t+1}
\ee
for each $n$ almost surely. 

By induction, we can conclude that for each $t\geq 0$, there exists $k^t$ such that $y_k(n)\geq M^t$ for all $k\geq k^t$ and each $n$. Since $M^t\rightarrow 0$ as $t\rightarrow \infty$, we obtain $\liminf_k y_k(n) \geq 0$ for each $n$.

\section{Proof of Lemma \ref{lem:bounds}}\label{app:bounds}

The proof follows from \eqref{eq:liminf}, which implies that given $(s,a)$, we have $\tilde{Q}_k^i(s,a) \geq (1-d_{\beta})^{-1}\underline{\eta}_k(s,a)$ for all $k$ for some asymptotically negligible term $\underline{\eta}_k(s,a)\rightarrow 0$ almost surely.
Particularly, we have
\be\label{eq:b1}
\bar{Q}_k = \tilde{Q}_k^1 + \tilde{Q}_k^2 \geq \tilde{Q}_k^1 + d_{\beta}\tilde{Q}_k^2 + \underline{\eta}_k = \Gamma_k + \underline{\eta}_k,
\ee
and
\be\label{eq:b2}
\Gamma_k = \tilde{Q}_k^1 + d_{\beta}\tilde{Q}_k^2 \geq d_{\beta}(\tilde{Q}_k^1+\tilde{Q}_k^2) - d_{\beta}\overline{\eta}_k = d_{\beta}\bar{Q}_k  - d_{\beta}\overline{\eta}_k
\ee
where we drop the arguments $(s,a)$ for notational convenience and $\overline{\eta}_k(s,a) = -\underline{\eta}_k(s,a)/d_{\beta}$. By \eqref{eq:b1} and \eqref{eq:b2}, we obtain
\be
\overline{\eta}_k + \frac{1}{d_{\beta}}\Gamma_k \geq \bar{Q}_k \geq \Gamma_k + \underline{\eta}_k.
\ee

\bibliography{mybibfile}

\begin{thebibliography}{28}
\providecommand{\natexlab}[1]{#1}
\providecommand{\url}[1]{\texttt{#1}}
\expandafter\ifx\csname urlstyle\endcsname\relax
  \providecommand{\doi}[1]{doi: #1}\else
  \providecommand{\doi}{doi: \begingroup \urlstyle{rm}\Url}\fi

\bibitem[Arslan and Yuksel(2017)]{ref:Arslan17}
G.~Arslan and S.~Yuksel.
\newblock Decentralized {Q}-learning for stochastic teams and games.
\newblock \emph{IEEE Transactions on Automatic Control}, 62\penalty0
  (4):\penalty0 1545--1558, 2017.

\bibitem[Benaim et~al.(2005)Benaim, Hofbauer, and Sorin]{ref:Benaim05}
M.~Benaim, J.~Hofbauer, and S.~Sorin.
\newblock Stochastic approximations and differential inclusions.
\newblock \emph{SIAM J. Control Optim.}, 44\penalty0 (1):\penalty0 328--348,
  2005.

\bibitem[Bogachev and Smolyanov(2020)]{ref:Bogachev20book}
V.~Bogachev and O.~G. Smolyanov.
\newblock \emph{Real and Functional Analysis}.
\newblock Springer Nature, 2020.

\bibitem[Borkar(2008)]{ref:Borkar08book}
V.~S. Borkar.
\newblock \emph{Stochastic Approximation: A Dynamical Systems Viewpoint}.
\newblock Hindustan Book Agency, 2008.

\bibitem[Chasnov et~al.(2020)Chasnov, Ratliff, Mazumdar, and
  Burden]{ref:Chasnov19}
B.~Chasnov, L.~Ratliff, E.~Mazumdar, and S.~Burden.
\newblock Convergence analysis of gradient-based learning in continuous games.
\newblock In \emph{Proceedings of The 35th Uncertainty in Artificial
  Intelligence Conference}, volume 115, pages 935--944, 2020.

\bibitem[Daskalakis et~al.(2020)Daskalakis, Foster, and
  Golowich]{ref:Daskalakis20}
C.~Daskalakis, D.~J. Foster, and N.~Golowich.
\newblock Independent policy gradient methods for competitive reinforcement
  learning.
\newblock In \emph{Advances in Neural Information Processing Systems}, 2020.

\bibitem[Filar and Vrieze(1997)]{ref:Filar97book}
J.~Filar and K.~Vrieze.
\newblock \emph{Competitive Markov Decision Processes}.
\newblock Springer Verlag, 1997.

\bibitem[Fink(1964)]{ref:Fink64}
A.~M. Fink.
\newblock Equilibrium in stochastic n-person game.
\newblock \emph{Journal of Science Hiroshima University Series A-I},
  28:\penalty0 89--93, 1964.

\bibitem[Fudenberg and Kreps(1995)]{ref:Fudenberg95}
D.~Fudenberg and D.~Kreps.
\newblock Learning in extensive-form games {I}. {S}elf-confirming equilibria.
\newblock \emph{Games and Economic Behavior}, 8:\penalty0 20--55, 1995.

\bibitem[Fudenberg and Levine(1998)]{ref:Fudenberg98}
D.~Fudenberg and D.~K. Levine.
\newblock \emph{The Theory of Learning in Games}.
\newblock MIT Press, Cambridge, MA, 1998.

\bibitem[Fudenberg and Levine(2009)]{ref:Fudenberg09}
D.~Fudenberg and D.~K. Levine.
\newblock Learning and equilibrium.
\newblock \emph{The Annual Review of Economics}, 1:\penalty0 385--419, 2009.

\bibitem[Harris(1998)]{ref:Harris98}
C.~Harris.
\newblock On the rate of convergence of continuous-time fictitious play.
\newblock \emph{Games and Economic Behavior}, 22:\penalty0 238--259, 1998.

\bibitem[Leslie and Collins(2003)]{ref:Leslie03}
D.~S. Leslie and E.~J. Collins.
\newblock Convergent multi-timescales reinforcement learning algorithms in
  normal form games.
\newblock \emph{The Annals of Applied Probability}, 13\penalty0 (4):\penalty0
  1231--1251, 2003.

\bibitem[Leslie and Collins(2005)]{ref:Leslie05}
D.~S. Leslie and E.~J. Collins.
\newblock Individual {Q}-learning in normal form games.
\newblock \emph{SIAM J. Control Optim.}, 44\penalty0 (2):\penalty0 495--514,
  2005.

\bibitem[Leslie et~al.(2020)Leslie, Perkins, and Xu]{ref:Leslie20}
D.~S. Leslie, S.~Perkins, and Z.~Xu.
\newblock Best-response dynamics in zero-sum stochastic games.
\newblock \emph{Journal of Economic Theory}, 189, 2020.

\bibitem[Littman(1994)]{ref:Littman94}
M.~L. Littman.
\newblock Markov games as a framework for multi-agent reinforcement learning.
\newblock In \emph{Proceedings of the 11th International Conference on Machine
  Learning (ICML)}, 1994.

\bibitem[Ozdaglar et~al.(2022)Ozdaglar, Sayin, and Zhang]{ref:OzdaglarICM}
A.~Ozdaglar, M.~O. Sayin, and K.~Zhang.
\newblock Independent learning in stochastic games.
\newblock In \emph{International Congress of Mathematicians}, 2022.

\bibitem[Sayin et~al.(2020)Sayin, Parise, and Ozdaglar]{ref:Sayin20}
M.~O. Sayin, F.~Parise, and A.~Ozdaglar.
\newblock Fictitious play in zero-sum stochastic games.
\newblock \emph{ArXiv:2010.04223}, 2020.

\bibitem[Sayin et~al.(2021)Sayin, Zhang, Leslie, Ozdaglar, and
  Ba\c{s}ar]{ref:Sayin21}
M.~O. Sayin, K.~Zhang, D.~Leslie, A.~Ozdaglar, and T.~Ba\c{s}ar.
\newblock Decentralized {Q}-learning in zero-sum markov games.
\newblock In \emph{Advances in Neural Information Processing Systems}, 2021.

\bibitem[Shapley(1953)]{ref:Shapley53}
L.~S. Shapley.
\newblock Stochastic games.
\newblock \emph{Proceedings of National Academy of Science USA}, 39\penalty0
  (10):\penalty0 1095--1100, 1953.

\bibitem[Sutton and Barto(2018)]{ref:Sutton18}
R.~S. Sutton and A.~G. Barto.
\newblock \emph{Reinforcement Learning: {A}n Introduction}.
\newblock MIT Press, Cambridge, MA, 2018.

\bibitem[Szepesvari and Littman(1999)]{ref:Szepesvari99}
C.~Szepesvari and M.~Littman.
\newblock A unified analysis of value-function-based reinforcement-learning
  algorithms.
\newblock \emph{Neural Computing}, 11:\penalty0 2017--2060, 1999.

\bibitem[Tsitsiklis(1994)]{ref:Tsitsiklis94}
J.~N. Tsitsiklis.
\newblock Asynchronous stochastic approximation and {Q}-learning.
\newblock \emph{Machine Learning}, 16:\penalty0 185--202, 1994.

\bibitem[Watkins and Dayan(1992)]{ref:Watkins92}
C.~J. C.~H. Watkins and P.~Dayan.
\newblock Q-learning.
\newblock \emph{Machine Learning}, 8\penalty0 (3):\penalty0 279--292, 1992.

\bibitem[Wei et~al.(2017)Wei, Hong, and Lu]{ref:Wei17}
{C.-Y.} Wei, {Y.-T.} Hong, and {C.-J.} Lu.
\newblock Online reinforcement learning in stochastic games.
\newblock In \emph{Proceedings of the 30th Conference on Neural Information
  Processing Systems (NIPS)}, 2017.

\bibitem[Young(2004)]{ref:Young04}
H.~P. Young.
\newblock \emph{Strategic Learning and Its Limits}.
\newblock Oxford University Press, 2004.

\bibitem[Zhang et~al.(2021)Zhang, Yang, and Ba\c{s}ar]{ref:Zhang19}
K.~Zhang, Z.~Yang, and T.~Ba\c{s}ar.
\newblock Multi-agent reinforcement learning: {A} selective overview of
  theories and algorithms.
\newblock In \emph{Handb. Rein. Learn. Cont.}, volume 325. Springer, 2021.

\bibitem[Zhu et~al.(2011)Zhu, Tembine, and Ba\c{s}ar]{ref:Zhu2011}
Q.~Zhu, H.~Tembine, and T.~Ba\c{s}ar.
\newblock Heterogeneous learning in zero-sum stochastic games with incomplete
  information.
\newblock In \emph{IEEE Conf. Decision and Control}, pages 219--224, 2011.

\end{thebibliography}

\end{document}